\documentclass{article}             
\usepackage{amssymb}
\usepackage[small,nohug,heads=vee]{diagrams}
\def\d{{\partial}}
\def\bV{{\mathbb{V}}}
\def\bR{{\mathbb{R}}}
\def\bS{{\mathbb{S}}}
\def\fg{{\mathfrak{g}}}
\def\fh{{\mathfrak{h}}}
\def\fm{{\mathfrak{m}}}

\newtheorem       {theorem}{Theorem}

\newtheorem{prop} [theorem]{Proposition}
\newtheorem{lemma}[theorem]{Lemma}

\begin{document}   

\begin{center}
{\Large{The affine approach to homogeneous geodesics\\
in homogeneous Finsler spaces}}
\bigskip

{\large{Zden\v ek Du\v sek}}
\bigskip

{\it{Dedicated to Professor Old\v rich Kowalski on the occasion of his 80th birthday}}
\end {center}
\bigskip

\begin{abstract} 
In a recent paper, it was claimed that any homogeneous Finsler space of odd dimension admits
a homogeneous geodesic through any point.
For the proof, the algebraic method dealing with the reductive decomposition of the Lie algebra
of the isometry group was used. However, the proof contains a serious gap.
In the present paper, homogeneous geodesics in Finsler homogeneous spaces are studied using the affine method,
which was developed in earlier papers by the author.
The mentioned statement is proved correctly and it is further proved that any homogeneous Berwald space
or homogeneous reversible Finsler space admits a homogeneous geodesic through any point.
\end{abstract}
\bigskip

\noindent
{\bf MSClassification:} {53C22, 53C60, 53C30}

\section{Introduction}
Let $M$ be either a pseudo-Riemannian manifold $(M,g)$, or a Finsler space $(M,F)$, or an affine manifold $(M,\nabla)$.
If there is a connected Lie group $G$ which acts transitively on $M$ as a group of isometries
or of affine diffeomorphisms, then $M$ is called a {\it homogeneous manifold\/}. It can be naturally identified
with the {\it homogeneous space\/} $(G/H,g)$, where $H$ is the isotropy group of the origin $p\in M$.

A geodesic $\gamma(s)$ through the point $p$ is {\it homogeneous} if it is an orbit of a one-parameter group
of isometries. More explicitly, if $s$ is an affine parameter and $\gamma(s)$ is defined in an open interval $J$,
there exists a diffeomorphism $s=\varphi(t)$ between the real line and the open interval $J$ and
a nonzero vector $X\in\fg$ such that $\gamma(\varphi(t))={\rm exp}(tX)(p)$ for all $t\in\bR$.
The vector $X$ is called a {\it geodesic vector}.
The diffeomorphism $\varphi(t)$ may be nontrivial only for null geodesics in a properly pseudo-Riemannian manifold
or for geodesics in affine manifolds.

A homogeneous Riemannian manifold $(M,g)$ or a homogeneous Finsler space $(M,F)$ is always
a {\it reductive homogeneous space\/}:
We denote by $\fg$ and $\fh$ the Lie algebras of $G$ and $H$ respectively
and consider the adjoint representation ${\rm Ad}\colon H\times\fg\rightarrow\fg$ of $H$ on $\fg$.
There exists a {\it reductive decomposition} of the form $\fg=\fm+\fh$ where $\fm\subset\fg$
is a vector subspace such that ${\rm Ad}(H)(\fm)\subset\fm$.
For a fixed reductive decomposition $\fg=\fm+\fh$ there is the natural identification
of $\fm\subset\fg=T_eG$ with the tangent space $T_pM$ via the projection $\pi\colon G\rightarrow G/H=M$.
Using this natural identification and the scalar product or the Finsler metric on $T_pM$, we obtain the invariant
scalar product $\langle\, ,\rangle $ or the invariant Minkowski norm $F$ and its fundamental tensor $g$ on $\fm$.
In the pseudo-Riemannian reductive case, geodesic vectors are characterized by the following {\it geodesic lemma}:
\begin{lemma}[\cite{KVa}, \cite{FMP}, \cite{DKa}]
 \label{golema}
Let $(G/H,g)$ be a reductive homogeneous pseudo-Riemannian manifold and $X\in{\fg}$.
Then the curve $\gamma(t)={\rm exp}(tX)(p)$ is geodesic with respect to some parameter $s$ if and only if
\begin{equation}
\nonumber
\langle [X,Z]_{\mathfrak m},X_{\mathfrak m}\rangle  = 
k\langle X_{\fm},Z\rangle
\end{equation}
for all $Z\in{\mathfrak m}$ and for some constant $k\in{\mathbb{R}}$.
If $k=0$, then $t$ is an affine parameter for this geodesic.
If $k\neq 0$, then $s=e^{-kt}$ is an affine parameter for the geodesic.
The second case can occur only if the curve $\gamma(t)$ is a null curve in a properly pseudo-Riemannian space.
\end{lemma}
The Finslerian version of this lemma was proved in \cite{La}:
\begin{lemma}[\cite{La}]
 \label{golema2}
Let $(G/H,g)$ be a homogeneous Finsler space. The vector $X\in{\fg}$ is a geodesic vector if and only if it holds
\begin{equation}
\nonumber
g_{X_\fm} ( [X,Z]_{\mathfrak m},X_{\mathfrak m} ) = 0
\end{equation}
for all $Z\in{\mathfrak m}$.
\end{lemma}

Another possible approach is to study the manifold $M$ using a more fundamental affine method,
which was proposed in \cite{D3} and \cite{DKV}.
It is based on the well known fact that a homogeneous manifold $M$ with the origin $p$
admits $n={\mathrm{dim}}M$ fundamental vector fields (Killing vector fields) which are linearly independent
at each point of some neighbourhood of $p$. Recall that a parametrized curve in a manifold $M$ is
{\it regular} if $\gamma\,'(t)\ne 0$ for all values of $t$.
It is well known that, in a homogeneous space $M=G/H$ with an invariant affine connection $\nabla$, each regular
orbit of a $1$-parameter subgroup $g_t\subset G$ on $M$ is an~integral curve of an~affine Killing vector field on $M$.
\begin{lemma}[\cite{DKV}]
\label{p:2} 
The integral curve $\gamma$ of a nonvanishing Killing vector field $Z$ on $M=(G/H,\nabla)$
is geodesic if and only if
\begin{eqnarray}
\label{mf1}
\nabla_{Z_{\gamma(t)}}Z=k_{\gamma}\cdot Z_{\gamma(t)}
\end{eqnarray}
holds along $\gamma$, where $k_{\gamma}\in\bR$ is a constant.
If $k_{\gamma}= 0$, then $t$ is the affine parameter of geodesic $\gamma$.
If $k_{\gamma}\neq 0$, then the affine parameter of this geodesic is $s=e^{k_\gamma t}$.
\end{lemma}

In the paper \cite{KS}, it was proved that any homogeneous Riemannian manifold admits a homogeneous geodesic
through the origin. The generalization to the pseudo-Riemannian (reductive and nonreductive) case was obtained
in \cite{D1} in the framework of a more general result, which says that any homogeneous affine manifold
$(M,\nabla)$ admits a homogeneous geodesic through the origin.
Here the affine method from \cite{DKV} and \cite{D3}, based on the study of integral curves of Killing vector fields,
was used. The proof is also using differential topology, namely smooth mappings $\bS^{n}\rightarrow \bS^n$.

In pseudo-Riemannian geometry, null homogeneous geodesics are of particular interest, see \cite{FMP} for instance.
In \cite{CM}, an example of a 3-dimensional Lie group with an invariant Lorentzian
metric which does not admit a light-like homogeneous geodesic was described.
Here the standard geodesic lemma was used, because the example is reductive.
In the paper \cite{D2}, the affine method used in \cite{D1}, \cite{D3} and \cite{DKV}
for the study of homogeneous affine manifolds was adapted to the pseudo-Riemannian case
and it was shown that any Lorentzian homogeneous manifold of even dimension admits a light-like
homogeneous geodesic through the origin.

Recently, in the paper \cite{ZY}, the existence of a homogeneous geodesic in homogeneous Finsler space of odd
dimension was claimed. The algebraic method developed in \cite{KS} and based on the reductive decomposition was
generalized to the Finserian situation and also differential topology and mappings $\bS^{n}\rightarrow \bS^n$ were used.
Surprisingly, neither the affine method nor the affine result from \cite{D1} was referred.
Moreover, the proof contains a serious gap.

In the present paper, the original result is reproved and the gap in the proof is indicated. It is shown
how the affine method can be adapted to the Finslerian setting and the mentioned result is proved correctly.
Further, it is proved that in homogeneous Berwald spaces and in homogeneous reversible Finsler spaces
a homogeneous geodesic always exists.

\section{Basic settings}
Recall that a {\it Minkowski norm} on the vector space $\bV$ is a nonnegative function $F:\bV\rightarrow\bR$
which is smooth on $\bV\setminus\{ 0 \}$, positively homogeneous ($F(\lambda y)=\lambda F(y)$ for any $\lambda >0$)
and whose Hessian $g_{ij} = (\frac{1}{2} F^2 )_{y^iy^j}$ is positively definite on $\bV\setminus \{ 0\}$.
Here $(y^i)$ are the components of a vector $y\in\bV$ with respect to a fixed basis $B$ of $\bV$\
and putting $y^i$ to a subscript means the patrial derivative.
Then the pair $(\bV,F)$ is called the {\it Minkowski space}.
The tensor $g_y$ with components $g_{ij}(y)$ is the {\it fundamental tensor}.
The {\it Cartan tensor} $C_y$ has components $C_{ijk}(y)=(\frac{1}{4} F^2)_{y^i y^j y^k}$.
A Finsler metric on the smooth manifold $M$ is a function $F$ on $TM$ which is
smooth on $TM\setminus \{ 0\}$ and whose restriction to any $T_xM$ is a Minkowski norm.
Then the pair $(M,F)$ is called the {\it Finsler space}.
On a Finsler space, functions $g_{ij}$ and $C_{ijk}$ depend smoothly on $x\in M$ and on $o\neq y\in T_xM$.

Further, we recall that the slit tangent bundle $TM_0$ is defined as $TM_0=TM\setminus\{ 0\}$.
Using the restriction of the natural projection $\pi\colon TM\rightarrow M$ to $TM_0$,
we naturally construct the pullback vector bundle $\pi^*TM$ over $TM_0$, as indicated in the following diagram:
 \begin{diagram}
  \pi^*TM  &            & TM \\
  \dTo     &            & \dTo{\pi} \\
  TM_0     & \rTo^{\pi} & M .
 \end{diagram}
For a given local coordinate system $(x^1,\dots,x^n)$ on $U\subset M$, at any $x\in M$, one has a natural basis
$\{\frac{\d}{\d x^1},\dots, \frac{\d}{\d x^n}\}$ of $T_xM$.
It is natural to express tangent vectors $y\in T_xM$ with respect to this basis.
Then $(x^i,y^i)$ is the {\it natural coordinate system} on $TU_0$.
We define further functions on $TU_0$, namely the {\it formal Christoffel symbols} $\gamma^i_{jk}$
and the {\it nonlinear connection} $N^i_j$, by the formulas
\begin{eqnarray}
\label{for4}
\gamma^i_{jk} & = & g^{is} \frac{1}{2} \bigl (
\frac{\d g_{sj}}{\d x^k} - \frac{\d g_{jk}}{\d x^s} + \frac{\d g_{ks}}{\d x^j} \bigr ), \cr
N^i_j & = & \gamma^i_{jk} y^k - C^i_{jk} \gamma^k_{rs} y^r y^s.
\end{eqnarray}
The {\it Chern connection} is the unique linear connection on the vector bundle
$\pi^*TM$ which is torsion free and almost $g$-compatible, hence its connection forms satisfy
\begin{eqnarray}
dx^j\wedge \omega^i_j & = & 0, \cr
dg_{ij} - g_{kj}\omega^k_i - g_{ik}\omega^k_j & = & 2C_{ijs}(dy^s + N^s_k dx^k).
\end{eqnarray}
It follows that it holds
\begin{eqnarray}
\label{for6}
\omega^i_j & = & \Gamma^i_{jk} dx^k, \cr
\Gamma^i_{jk} & = & \Gamma^i_{kj}, \cr
\Gamma^l_{jk} & = & \gamma^l_{jk} - g^{li} ( C_{ijs} N^s_k - C_{jks} N^s_i + C_{kis} N^s_j),
\end{eqnarray}
see some monograph, for example \cite{BCS} or \cite{De} for details.
If we fix a nowhere vanishing vector field $V$ on $M$, we obtain an affine connection $\nabla^V$ on $M$.
In the local chart, it is expressed with respect to arbitrary vector fields
$W_1=W_1^i\frac{\d}{\d x^i}$ and $W_2=W_2^i\frac{\d}{\d x^i}$ by the formula
\begin{eqnarray}
\label{f1}
\nabla^V_{W_1} W_2 |_x = \bigl [ {W_1}(W_2^i) + W_2^j {W_1}^k \Gamma^i_{jk}(x,V) \bigr ] \frac{\d}{\d x^i}.
\end{eqnarray}
The affine connection $\nabla^V$ on $M$ is torsion free and almost metric compatible, which means
\begin{eqnarray}
\label{perp}
\nabla^V_{W_1} {W_2} - \nabla^V_{W_2} {W_1} & = & [W_1,W_2],\cr
W g_V(W_1,W_2) & = & g_V(\nabla^V_W W_1,W_2) + g_V(W_1, \nabla^V_W W_2) + \cr
                  && + 2C_V (\nabla^V_W V, W_1, W_2) ,
\end{eqnarray}
for arbitrary vector fields $W, W_1, W_2$.
Using the affine connection $\nabla^V$, we define the derivative along a curve $\gamma(t)$
with velocity vector field $T$. Let $W_1,W_2$ be vector fields along $\gamma$, we define
\begin{eqnarray}
\label{f3}
D_{W_1}W_2 = \nabla^{T'}_{W_1'} {W_2'},
\end{eqnarray}
where the vector fields $T'$, $W_1'$ and $W_2'$ on the right-hand side are smooth extensions of $T$, $W_1$ and $W_2$
to the neighbourhood of $\gamma(t)$. The definition above does not depend on the particular extension.
A regular smooth curve $\gamma$ with tangent vector field $T$ is a {\it geodesic} if $D_{T} (\frac{T}{F(T)}) = 0$.
In particular, a geodesic of constant speed satisfies $D_{T} {T} = 0$.

We now reprove with minor technical modifications the result from \cite{ZY}.
It is using the algebraic method developed in \cite{KS} for the Riemannian metric.
\begin{theorem}[\cite{ZY}]
\label{ty}
Let $(M,F)$ be a homogeneous Finsler space of odd dimension and $p\in M$.
Then $M$ admits a homogeneous geodesic through $p$.
\end{theorem}
{\it Proof.}
Let $G$ be a group of isometries acting transitively on $M$ and $H$ the isotropy subgroup of the origin $p\in M$.
We can write $M=G/H$. Denote by ${\mathfrak{g}}$ and ${\mathfrak{h}}$ the corresponding Lie algebras,
by $K$ the Killing form of ${\mathfrak{g}}$ and by ${\mathrm{rad}}(K)$ the null space of $K$.
Because $K$ is nondegenerate on ${\mathfrak{h}}$, we can put ${\mathfrak{m}}={\mathfrak{h}}^\perp$
and fix the reductive decomposition ${\mathfrak{g}}= {\mathfrak{h}}+{\mathfrak{m}}$.
It holds ${\mathrm{rad}}(K)\subseteq {\mathfrak{m}}$ and there are the two possible cases:

If ${\mathrm{rad}}(K)={\mathfrak{m}}$, then
$[{\mathfrak{g}}, {\mathfrak{g}}]_{\mathfrak{m}}$ is a proper subset of ${\mathfrak{m}}$.
We choose arbitrary vector $X\in [{\mathfrak{g}}, {\mathfrak{g}}]_{\mathfrak{m}}^\perp$ and
for any $Z\in\fg$ it holds $[X,Z]_\fm \in [\fg,\fg]_\fm$, hence
\begin{eqnarray}
\nonumber
g_X(X,[X,Z]_{\mathfrak{m}}) =0\quad\forall Z\in{\mathfrak{g}}
\end{eqnarray}
and $X$ is a geodesic vector.

If ${\mathrm{rad}}(K)\subsetneq {\mathfrak{m}}$, we fix an invariant scalar product $\langle, \rangle$ on $\fm$.
For each unit vector $X\in {\mathrm{S}}^{n-1}\subset\fm$, we define the operator
$\alpha^X\colon\fm\rightarrow\fm$ by the formula
\begin{eqnarray}
\nonumber
g_X(\alpha^X U,V) = K(U,V) \quad \forall U,V\in\fm.
\end{eqnarray}
If the Finsler space is Riemannian, there is just one operator $\alpha$ and we can continue as in \cite{KS}:
There always exists a nonzero eigenvector $Y$ of $\alpha$ with nonzero eigenvalue $\lambda$,
because ${\mathrm{rad}}(K)\subsetneq {\mathfrak{m}}$. We obtain
\begin{eqnarray}
\label{kv}
g(Y, [Y,Z]_{\mathfrak{m}}) & = & \frac{1}{\lambda} g(\alpha(Y), [Y,Z]_{\mathfrak{m}}) = 
\frac{1}{\lambda} K(Y, [Y,Z]_{\mathfrak{m}}) = \cr
& = & \frac{1}{\lambda} K(Y, [Y,Z]) = \frac{1}{\lambda} K([Y ,Y],Z]) = 0 \quad \forall Z\in\fm
\end{eqnarray}
and $Y$ is a geodesic vector. We have used here also the invariance of the Killing form $K$.
In a general Finsler space, if there is a vector $\bar X\in \fm$ such that the eigenvector
$Y^{\bar X}$ of $\alpha^{\bar X}$ satisfies $Y^{\bar X}=\bar X$, we can use similar
steps as in the formula (\ref{kv}) above, write
\begin{eqnarray}
\nonumber
g_{\bar X}(\bar X, [\bar X,Z]_{\mathfrak{m}}) 
& = & \frac{1}{\lambda(\bar X)} g_{\bar X}(\alpha^{\bar X}(\bar X), [\bar X,Z]_{\mathfrak{m}})
  =   \frac{1}{\lambda(\bar X)} K(\bar X, [\bar X,Z]_{\mathfrak{m}}) = \cr
& = & \frac{1}{\lambda(\bar X)} K(\bar X, [\bar X,Z]) = \frac{1}{\lambda(\bar X)} K([\bar X ,\bar X],Z])
  = 0 \quad \forall Z\in\fm
\end{eqnarray}
and $\bar X$ is a geodesic vector. In the paper \cite{ZY}, the mapping
$v\colon {\mathrm{S}}^{n-1}\rightarrow {\mathrm{S}}^{n-1}$ was constructed by the assignment $X\mapsto Y^X$,
where $Y^X$ is the eigenvector of the operator $\alpha^X$ with maximal absolute value of the eigenvalue $\lambda(X)$.
The mapping $v$ was claimed to be continuous and the fixed point theorem was used.
However, this part of the proof is not well justified and probably it is wrong.
For the family of operators $\alpha^X$, the assignment of the eigenvector $Y^X$ with maximal absolute value
of the eigenvalue $\lambda(X)$ is not a continuous mapping in general.
This seems to be a serious gap in the proof, because there is not an obvious way how to correct it.
$~\hfill\square$

\section{Affine method for Finsler spaces}
First, let us formulate simple observations which follow from homogeneity of the Finsler metric $F$.
\begin{prop}
\label{p1}
Let $(M,F)$ be a homogeneous Finsler space, $G$ a group of isometries acting transitively on $M$, $X^*$ a Killing
vector field generated by the vector $X\in\fg$, $\phi(t)={\mathrm{exp}}(tX)$ and $\gamma(t)$ the integral curve
of $X^*$ through $p\in M$. Along the curve $\gamma(t)$, it holds
\begin{eqnarray}
\phi(t)(p) & = & \gamma(t), \cr
\phi(t)_*(X^*(p)) & = & X^*(\gamma(t))
\end{eqnarray}
and
\begin{eqnarray}
F(\phi(t)(p),\phi(t)_*V) & = & F(p,V), \cr
g_{(\gamma(t),X^*(\gamma(t)))}(\phi(t)_*U,\phi(t)_*V) & = & g_{(p,X^*(p))}(U,V),
\end{eqnarray}
for all $t\in\bR$ and for all $U, V \in T_pM$.
\end{prop}
\begin{prop}
\label{p6}
With the same assumptions as in Proposition \ref{p1}, along the curve $\gamma(t)$, it holds
\begin{eqnarray}
g_{(\gamma(t),X^*(\gamma(t))}(D_{X^*}X^*\big |_{\gamma(t)},\phi(t)_*U) & = & g_{(p,X^*(p))}(D_{X^*}X^*\big |_p,U),
\end{eqnarray}
for all $t\in\bR$ and for all $U \in T_pM$. Consequently, if
\begin{eqnarray}
\label{f6}
D_{X^*}X^* \big |_p= 0,
\end{eqnarray}
then the curve $\gamma(t)$ is a homogeneous geodesic.
\end{prop}
We shall now give a correct proof of Theorem \ref{ty}.
\begin{theorem}
\label{t7}
Let $(M,F)$ be a homogeneous Finsler space of odd dimension and $p\in M$.
Then $M$ admits a homogeneous geodesic through $p$.
\end{theorem}
{\it Proof.}
Let us consider the Killing vector fields $K_1,\dots,K_n$ which are linearly independent at each point
of some neighbourhood $\mathcal U$ of $p$ and denote by $B$ the basis $\{K_1(p),\dots,K_n(p)\}$ of $T_pM$.
Any tangent vector $X\in T_pM$ has coordinates $(x_1,\dots x_n)$ with respect to the basis $B$.
These coordinates determine the Killing vector field $X^*=x_1K_1+\dots+x_nK_n$ and an integral curve
$\gamma$ of $X^*$ through $p$. We are going to show that there exists a vector $\bar X\in T_pM$ such that
the integral curve $\gamma$ of $\bar X^*$ through $p$ is geodesic.

Let us consider the sphere $\bS^{n-1}$ of vectors $X\in T_pM$ whose coordinates $(x_1,\dots,x_n)$ with respect to $B$
have the norm equal to 1 with respect to the standard Euclidean scalar product $\langle , \rangle$ on $\bR^n$.
In other words, the scalar product $\langle, \rangle$ is chosen in a way that the above basis $B$ is orthonormal.
We stress that this scalar product does not come from any Finslerian product $g$ used so far.
For each $X\in \bS^{n-1}$, denote by $v(X)$ the derivative $D_{X^*_{\gamma(t)}}X^*|_{t=0}$.
Further, denote by $t(X)$ the vector $v(X)-\langle v(X),X\rangle X$.
Then, for each $X\in \bS^{n-1}$, $t(X)\perp X$ with respect to the above Euclidean scalar product.
Clearly, the map $X\mapsto t(X)$ defines a smooth tangent vector field on the sphere $\bS^{n-1}$.
If $n$ is odd, according to a well known fact from differential topology, there is a vector $\bar X$ such that $t(\bar X)=0$.

To finish the proof, we use the formula (\ref{perp}) and the standard fact that $C_{X^*}(X^*,X^*,X^*)=0$.
We observe that, for each $X\in \bS^{n-1}\subset T_pM$, it holds
\begin{eqnarray}
g_{(p,X)}(v(X),X) = g_{(p,X)} (D_{X^*_{\gamma(t)}}X^*\big |_{t=0},X^*_p) & = & 0
\end{eqnarray}
and hence $v(X)$ lies in the orthogonal complement of $X$ in $T_pM$ with respect to the scalar product $g_{(p,X)}$.
The vector $t(X)$ is the projection of $v(X)$ to another complementary subspace of $X$ in $T_pM$ and hence $v(X)=0$
if and only if $t(X)=0$. If follows, using also Proposition \ref{p6} and formula (\ref{f6}), that the integral curve
of the vector field $\bar X^*$ through $p$ is a homogeneous geodesic.
$~\hfill\square$
\medskip

Let us now recall that the Finsler metric $F$ is called a {\it Berwald metric} if the Christoffel symbols
$\Gamma^i_{jk}(x,y)$ of the Chern connection in natural coordinates do not depend on the direction $y$, hence
$\Gamma^i_{jk}(x,y)= \Gamma^i_{jk}(x)$.
We further recall that the Finsler metric $F$ is {\it reversible} if, for any point $x\in M$ and for any vector $y\in T_xM$,
it holds $F(x,y)=F(x,-y)$.
\begin{theorem}
Let $(M,F)$ be a homogeneous Berwald space or a homogeneous reversible Finsler space and let $p\in M$.
Then $M$ admits a homogeneous geodesic through $p$.
\end{theorem}
{\it Proof.}
If the Finsler metric $F$ is Berwald, using formulas (\ref{f1}) and (\ref{f3})
we easily deduce that for any Killing vector field $X^*$ it holds
\begin{eqnarray}
\label{main}
D_{X^*} {X^*}  = \nabla^{X^*}_{X^*} {X^*} = \nabla^{-X^*}_{-X^*} {-X^*} = D_{-X^*} {-X^*}.
\end{eqnarray}
For a reversible Finsler metric, one can check by the straightforward calculations
and using formula (\ref{for4}), that it holds
\begin{eqnarray}
g_{ij}(x,y)  =  g_{ij}(x,-y), && C_{ijk}(x,y) =  - C_{ijk}(x,-y), \cr
\gamma^i_{jk}(x,y) = \gamma^i_{jk}(x,-y), && N^i_{j}(x,y) = - N^i_{j}(x,-y).
\end{eqnarray}
Further, using formula (\ref{for6}), we obtain
\begin{eqnarray}
\Gamma^i_{jk}(x,y) & = & \Gamma^i_{jk}(x,-y)
\end{eqnarray}
and using formula (\ref{f1}), we obtain again that formula (\ref{main}) is valid also in this situation.
Formula (\ref{main}) is essential for the next step.

Let us use the same notation and setting as in the proof of Theorem \ref{t7} and let us consider the mappings
$v\colon \bS^{n-1}\rightarrow T_pM$ and $t\colon\bS^{n-1}\rightarrow T_pM$.
For $n$ even, let us assume that $t(X)\neq 0$ everywhere.
Putting $f(X)=t(X)/\|t(X)\|$, where the norm comes from the Euclidean scalar product $\langle,\rangle$,
we obtain a smooth map $f\colon \bS^{n-1}\rightarrow \bS^{n-1}$ without fixed points.
According to a well known statement from differential topology, the degree of $f$ is $\deg(f)=(-1)^n$,
because it is homotopic to the antipodal mapping.
On the other hand, according to formula (\ref{main}), we have
$v(X)=v(-X)$ and hence $f(X)=f(-X)$ for each $X$. If $Y$ is a regular value of $f$, then the inverse image
$f^{-1}(Y)$ consists of even number of elements. Hence $\deg(f)$ is an even number, which is a contradiction.
Hence, the assumption $t(X)\neq 0$ was wrong.

It follows that there is again a vector $\bar X\in T_pM$ such that $t(\bar X)=0$ and also $v(\bar X)=0$.
Consequently, the integral curve of the vector field $\bar X^*$ through $p$ is a homogeneous geodesic.
$~\hfill\square$

% \begin{center}
% {\bf Acknowledgements}
% \end{center}
% The authors were supported by the grant GA\v{C}R 14-02476S.
% \bigskip

\noindent
Address of the author:\\
Zden\v ek Du\v sek\\
University of Hradec Kr\'alov\'e, Faculty of Science\\
Rokitansk\'eho 62, 500 03 Hradec Kr\'alov\'e, Czech Republic\\
zdenek.dusek@uhk.cz

\end{document}